\documentclass[12pt]{article}
\usepackage{amssymb}
\usepackage{color}
\usepackage{graphicx}

\newcommand{\C}{\mathbb C}
\newcommand{\D}{\mathbb D}
\newcommand{\F}{\mathbb F}
\newcommand{\N}{\mathbb N}
\newcommand{\R}{\mathbb R}
\newcommand{\T}{\mathbb T}
\newcommand{\Rp}{{\mathbb R}_+^\times}
\newcommand{\cS}{\mathcal X}

\renewcommand{\H}{\mathbb H}
\renewcommand{\P}{\mathbb P}
\newcommand{\Q}{\mathbb Q}
\newcommand{\Z}{\mathbb Z}

\newcommand{\re}{\mathrm{Re}}
\newcommand{\im}{\mathrm{Im}}

\newcommand{\ex}{\mathbf{e}}
\renewcommand{\a}{\alpha}

\newcommand{\g}{\gamma}

\newcommand{\De}{\mathit{\Delta}}

\newcommand{\h}{\vartheta}
\newcommand{\hh}{\mathit{\Theta}}
\newcommand{\e}{\varepsilon}
\newcommand{\f}{\varphi}
\newcommand{\z}{\zeta}
\newcommand{\G}{\mathit{\Gamma}}
\newcommand{\bG}{\mathbf{\Gamma}}

\newcommand{\w}{\omega}
\renewcommand{\i}{\mathbf{i}}
\newcommand{\Mp}{M^\times}
\newcommand{\la}{\langle}
\newcommand{\ra}{\rangle}
\newcommand{\tr}{\;^t}
\newcommand{\diag}{\mathrm{diag}}

\newcommand{\iint}{\int\!\!\!\!\int}
\newcommand{\pr}{\textit{Proof.}\hspace{2mm} }
\newcommand{\qed}{\hfill{$\square$} \bigbreak}
\newcommand{\KK}{\textsl K3 }
\newtheorem{theorem}{Theorem}
\newtheorem{proposition}{Proposition}
\newtheorem{lemma}{Lemma}

\newtheorem{cor}{Corollary}

\newtheorem{remark}{Remark}

\def\comment#1{}
\def\ul#1{\underline{#1}}
\title
{Thomae type formula for \KK surfaces given by 
double covers of the projective plane 
branching along six lines}
\author{Keiji \textsc{Matsumoto} and
Tomohide \textsc{Terasoma}
}
\begin{document}
\maketitle   
\begin{center}
\textit{
Dedicated to Professor Takayuki Oda on his sixtieth birthday
}
\end{center}
\begin{abstract}
In this paper, we give Thomae type formula for \KK surfaces $\cS$ 
given by double covers of the projective plane branching along six lines. 
This formula gives relations 
between theta constants on the bounded symmetric domain 
of type $I_{22}$ and period integrals of 
$\cS$.
Moreover, we express the period integrals by using 
the hypergeometric function $F_S$  of four variables.
As applications of our main theorem, we define $\R^4$-valued sequences 
by mean iterations of four terms, and express their common limits 
by the hypergeometric function $F_S$. 

\end{abstract}
\setcounter{footnote}{1}


\maketitle
\tableofcontents
\section{Introduction}
Let us consider period integrals
\begin{equation}
\label{eq:definition of elliptic period}
\w_A(\lambda)=\int^{1/\lambda}_1\frac{dt}{\sqrt{t(1-t)(1-\lambda t)}}
,\quad
\w_B(\lambda)=\int_{0}^1\frac{dt}{\sqrt{t(1-t)(1-\lambda t)}},
\end{equation}
of an elliptic curve $s^2=t(1-t)(1-\lambda t)$ 
with  $\lambda\in\C-\{0,1\}$. 
If $\lambda$ belongs to 
the open interval $(0,1)$, then 
they are expressed by  
the Gauss hypergeometric function  $F(a,b,c;z)= 
\sum\limits_{n=0}^\infty{(a)_n(b)_n\over (c)_nn!}z^n:$
$$
\w_A(\lambda)=\i\pi F({1\over 2},{1\over 2},1;1-\lambda),
\quad 
\w_B(\lambda)=\pi F({1\over 2},{1\over 2},1; \lambda),
$$
where $\i=\sqrt{-1}$. 
%
The function $\tau=\w_A(\lambda)/\w_B(\lambda)$ of $\lambda$ is
continued to a map
$$
per: \widetilde{X}\to 
\H=\{z \in \C\mid \im(z)>0\}
$$
from the universal covering $\widetilde{X}$ of $\C-\{0,1\}$
to $\H$, which is called the period map.
The inverse of the map $per$ can be described as
$$
\lambda={\h_{[10]}^4(\tau)/\h_{[00]}^4(\tau)},
$$ 
where
$$
\h_{[ab]}(\tau)
=\sum_{n\in \Z}\exp[\pi\i \{(n+{a\over2})^2\tau+ (n+{a\over2})b\}]\quad 
([ab]=[00],[01],[10])
$$
is Jacobi's theta constant.
Under these correspondences of variables $\lambda\in \C-\{0,1\}$ and 
$\tau\in \H$,  
the theta constant and the elliptic integral are related as
\begin{equation}
\label{Jacobi's formula}
\vartheta_{[ab]}^4(\tau)=
\frac{\Lambda_{[ab]}}{\pi^2}\w_B(\lambda)^2,
\end{equation}
where 
$$
\Lambda_{[00]}=1, \quad \Lambda_{[01]}=1-\lambda,\quad 
\Lambda_{[10]}=\lambda.
$$
The identity (\ref{Jacobi's formula}) is called Jacobi's formula.
On the other hand, we have
the $2\tau$-formulas for the theta constants
$$
\h_{[00]}^2(2\tau)=\frac{\h_{[00]}^2(\tau)+\h_{[01]}^2(\tau)}{2},
\quad 
\h_{[01]}^2(2\tau)=\h_{[00]}(\tau)\h_{[01]}(\tau).
$$

These formulas are applied to the study of the 
arithmetic-geometric mean
as follows.
Let $c_1,c_2\in \Rp$ be positive real numbers. We define vector valued sequence
$\{m^n(c_1,c_2)\}_{n\in \N}$
by 
$$
m^n(c_1,c_2)=
\overbrace{m\circ \cdots \circ m}^n(c_1,c_2)
$$
where the map $m:(\Rp)^2 \to (\Rp)^2$ is
\begin{eqnarray*}
m(u_1,u_2)=
(\displaystyle \frac{u_1+u_2}{2}, 
\sqrt{u_1 u_2}).
\end{eqnarray*}
Both components have a common limit and it is called the
arithmetic-geometric mean and denoted by
$m^\infty_*(c_1,c_2)$.
Using Jacobi's formula and $2\tau$-formulas, we have a relation between
the arithmetic-geometric mean and the hypergeometric function:
$$m^\infty_*(c_1,c_2)=
{c_1\over F({1\over 2},{1\over 2},1;1-({c_2\over c_1})^2)}.$$
By the above relation and the invariance property 
$m^\infty_*(m(c_1,c_2)))=m^\infty_*(c_1,c_2)$, we have
the Gauss transformation formula
\begin{equation}
\label{gauss transform classical}
F\Big({1\over 2},{1\over 2},1;1-{4z\over (1+z)^2}\Big)
={1+z\over 2}F\big({1\over 2},{1\over 2},1;1-z^2\big).
\end{equation}

Thomae studies period integrals of a hyperelliptic curve of arbitrary genus  
and generalizes Jacobi's formula to Thomae's formulas in \cite{To}.
Based on $2\tau$-formulas of theta constants defined on the Siegel upper half
space $\H_2$ of degree $2$, 
Borchardt introduces a vector valued sequence $\{m^n(c_1,\dots,c_4)\}_{n\in \N}$ 
with initial $(c_1,\dots,c_4)\in (\Rp)^4$ 
given by the iteration of the map 
$$m:(\Rp)^4\ni u=(u_1,\dots,u_4)\mapsto 
(m_1(u),\dots,m_4(u))\in (\Rp)^4,$$
where 
\begin{eqnarray*}
m_1(u)&=&{u_1+u_2+u_3+u_4\over 4},\quad
m_2(u)={\sqrt{u_1u_2}+\sqrt{u_3u_4}\over 2},\\
m_3(u)&=&{\sqrt{u_1u_3}+\sqrt{u_2u_4}\over 2},\hspace{8mm}
m_4(u)={\sqrt{u_1u_4}+\sqrt{u_2u_3}\over 2}.
\end{eqnarray*}
By using Thomae's formulas,  
he expresses the common limit
of the components of
$\{m^n(c_1,\dots,c_4)\}_{n\in \N}$
by period integrals of a genus $2$ hyperelliptic curve.
For related studies,
refer to \cite{B}, \cite{MT} and \cite{Me}.

In this paper, we give Thomae type formula for
\KK surfaces which are double covers of 
the complex projective plane $\P^2$ branching along normal crossing six lines.
The configurations of normal crossing six lines are parametrized 
by $3\times 6$ matrices
$x$ and the corresponding \KK surface is denoted by $\cS(x)$.
Period integrals of $\cS(x)$ are expressed
in terms of two kinds of hypergeometric functions $F_S$ and $F_T$ 
of four variables defined in (\ref{S-series}) and (\ref{T-series}), respectively.
In \S \ref{isom symmetric domain}, 
we define a normalized period matrix $\tau$ of $\cS(x)$ in
the $4$-dimensional bounded symmetric domain $\D$ 
of type $I_{22}$.
Let $P_{3,3}$ be the set of unorderd pair $\la J\ra =(J,J^c)$, such that
$\# J=\#J^c=3$ and $J\cup J^c=\{1,\dots,6\}$. Then we have
$\# P_{3,3}=10$.
To state the main theorem,
we introduce the following notations:
\begin{enumerate}
\item
$\hh_{\la J \ra}(\tau)$ is
theta functions 
on $\D$
evaluated at the normalized 
period matrix $\tau$ indexed by $\la J\ra \in P_{3,3}$,
\item
$x\la J \ra$ is the product of two $3\times3$-minors 
of a $3\times 6$-matrix $x$ in the configuration space $X(3,6)$ also indexed 
by $\la J\ra \in P_{3,3}$,
\item
$\w_{34}(x)$
is the period integral 
of the \KK surface $\cS(x)$ given in (\ref{periods}).
\end{enumerate}
Then the main theorem is the identity 
\begin{equation}
\label{form of main Th}
\hh_{\la J \ra}^2(\tau)=
{1\over 4\pi^4}x\la J \ra \w_{34}(x)^2
\end{equation}
for any $\la J\ra \in P_{3,3}$.

The subfamily consisting of Kummer varieties of principally
polarized abelian varieties is called the Kummer locus.
This locus corresponds to the Siegle
upper half space $\H_2$ realized as a closed subdomain of $\D$.
Our identity becomes Thomae's
formula for genus two curves on this locus.

In the paper \cite{MSY} and \cite{Ma}, they prove that
the point $[\hh_{\la J \ra}(\tau)]_{\la J\ra \in P_{3,3}}$
in $\P^9$ is equal to
$[x\la J \ra]_{\la J\ra \in P_{3,3}}$.
The key for our proof of the main theorem is the study of the relation 
between a period of $\cS$ 
and the automorphic factor of $\hh_{\la J \ra}$ 
by the action of the monodromy group of $per$ 
via the isomorphism between $\D$ and 
$D_H$ defined in \S \ref{isom symmetric domain}.

As an application of our main theorem, we study 
a vector valued sequence 
obtained by mean iteration of a map from $\R^4_{+}$ to $\R^4_{+}$ 
which is different from that defined by Borchardt
in \S \ref{our AGM}.
We show that this vector valued sequence has a common limit 
and that it can be expressed by the hypergeometric function $F_S$.
This formula is obtained by the main theorem and 
$2\tau$ formulas for the theta functions $\hh_{\la J \ra}$ 
in Theorem \ref{2tau-formula}.
We also give an explanation on the relation between Borchardt's 
arithmetic-geometric mean $m^\infty_*(c_1,\dots,c_4)$ and
the hypergeometric function $F_S$ in 
\S \ref{Borchardt's AGM}.
In the last section, we prove 
several functional equations of the hypergeometric function $F_S$
arising from the invariance property for vector valued 
mean iterations.
These are analogs of the Gauss transformation formula 
(\ref{gauss transform classical}) 
for the hypergeometric function $F_S$.

\section{Certain family of K3 surfaces}
\subsection{Double coverings of $\P^2$ branching along 6 lines}
Let $\Mp(3,6)$ be the open subset of $M(3,6)$ defined by
$$
\Mp(3,6)=\Big\{x=
(\ell_1, \dots, \ell_6)
\in M(3,6) \Big|
\begin{minipage}{5cm}
\textrm{the determinants of $(3,3)$- } 
\textrm{minors are non-zero }
\end{minipage} \Big\}
$$
For $\ell_i=\!\!\ ^t(\ell_{0i},\ell_{1i},\ell_{2i})$,
we define a linear function $(t,\ell_i)$ by $\sum_{j=1}^3\ell_{j,i}t_j$.
Let $\widehat{\cS}^*$ be the variety defined by
$$
\widehat{\cS}^*=\{
(t:y)\times x
\in \P(1,1,1,3) \times
\Mp(3,6)
\mid
y^2=(t,\ell_1)\cdots (t,\ell_6)
\},
$$
where $(t:y)=(t_0:t_1:t_2:y),$ $x=(\ell_1, \dots, \ell_6)$ and 
$\P(1,1,1,3)$ is the weighted projective space of weight $(1,1,1,3)$.
Then by the natural map $pr_1:\widehat\cS^*\to \P^2 \times \Mp(3,6)$,
$\widehat\cS^*$ is a family of branched covering of $\P^2$ over $\Mp(3,6)$.
By resolving singularities, we have a family of K3 
surfaces $pr_2:\widehat\cS \to \Mp(3,6)$
on $\Mp(3,6)$. 

Let $\T$ be a torus defined by 
$$
\T=\{\lambda=(\lambda_0,\lambda_1, \dots, \lambda_6)\in {\C^{\times}}^6\times \C^{\times}\mid
\lambda_0^2=\lambda_1 \cdots \lambda_6\}.
$$
The group $GL_3(\C)\times \T$ acts on $\widehat{\cS}$ 
by 
$$
t\mapsto tg^{-1},\quad x\mapsto g\cdot x\cdot \diag(\lambda_1, \dots,
\lambda_6),
\quad y \mapsto \lambda_0 y
$$
for $(g,\lambda)\in GL_3(\C)\times \T$ and
it induces an action of $GL_3(\C)\times \T$ on $\Mp(3,6)$.
The quotients 
of $\widehat{\cS}$ and $\Mp(3,6)$ by $GL_3(\C)\times \T$
are denoted by $\cS$ and $X$, respectively.
The natural map
$pr_2:\widehat{\cS} \to \Mp(3,6)$ induces a map $\cS \to
X$, which is also denoted by $pr_2$.
The variety $X$ is equal to 
the double coset space: 
$$
X=X(3,6)=GL_3(\C)\backslash \Mp(3,6)/(\C^\times)^6,
$$
which is called the configuration space.
The fiber of $\cS$ at $x\in X$ is denoted by $\cS(x)$.

There are $15$ rational curves $\tilde l_{jk}(x)$ $(1\le j< k\le 6)$  
in $\cS(x)$ coming from the resolutions of nodes at 
$l_{jk}(x)=l_j(x)\cap l_k(x)$, where $l_i(x)$ is the line defined
by $(t,\ell_i)=0$.
Let $\tilde l_0(x)$ be a pull back of a generic line in $\P^2$ by $pr_1$.
Let $S(x)$ be the subgroup of $H_2(\cS(x),\Z)$ generated by 
algebraic cycles $\tilde l_{jk}(x)$ and $\tilde l_0(x)$. 
Its orthogonal complement $T(x)$ in $H_2(\cS(x),\Z)$ with respect to 
the intersection form $(\; \cdot\;  )$ is called the transcendental
lattice of $X$ and its rank is $22-16=6$. 
In \S \ref{topological cycle at ref pt}, 
we give a basis of $T(x)$ and its dual in $H_2(\cS(x),\Z)$. 
These bases are slightly different from those defined in \cite{MSY} and \cite{Y}.

\subsection{Relative invariants and a global two form}

The characters 
$$
\rho:GL_3(\C)\times \T\ni
(g, \lambda
)\mapsto \deg(g)\lambda_0
\in
\C^{\times}:
$$
and $\rho^2$
define linearlizations of 
$GL_3(\C)\times \T$ of $\mathcal O_{\Mp(3,6)}$ and
$\mathcal O_{\widehat \cS}$.
The invariant line bundles on $\cS$
and $X$ under these actions are denoted by 
$\mathcal L$ and $\mathcal M$, respectively. 
We have $\mathcal L^{\otimes 2}=pr_2^*\mathcal M$.

We construct elements of $H^0(X,\mathcal M)$.
Let $J$ be a subset of the set $\{1,\dots,6\}$ with cardinality $3$ 
and $J^c$ be its complement.
By reordering elements, we may write $J$ and $J^c$ as 
$$J=\{j_1,j_2,j_3\},\ j_1<j_2<j_3,\quad J^c=\{j_4,j_5,j_6\},\ j_4<j_5<j_6.$$
A pair $\la J\ra=(J,J^c)=(J^c,J)$ of $J$ and $J^c$ 
is called a $(3,3)$-partition of the set $\{1,\dots,6\}$. 
The set of $(3,3)$-partitions is denoted by $P_{3,3}$.
Note that $\#P_{3,3}=10$.  
For $x=(x_{ij})\in M(3,6)$  and $\la J\ra \in P_{3,3}$, 
we set 
$$
x\la J\ra=\det(x_{i,j_k})_{1\le i,k\le 3}\det(x_{i,j_{k+3}})_{1\le i,k\le 3}
.$$
Then $x\la J\ra$ is an element of $H^0(X, \mathcal M)$.
By Pl\"ucker relations, we have the following. 
\begin{lemma} \label{Pl-rel}
Let $St$ be the set of (2,2,2)-standard tableaux i.e.
$(J,J^c)=(\{j_1,j_2,j_3\},$
$\{j_4,j_5,j_6\})$ with
$$
\begin{array}{ccc}
j_1 &<& j_4 \\
\wedge & & \wedge \\ 
j_2 &<& j_5 \\ 
\wedge & & \wedge \\ 
j_3 &<& j_6. \\ 
\end{array}
$$
Then $\# St=5$ and $\{x\la J \ra\mid \la J\ra \in St\}$ forms a basis
of a linear system 
in $H^0(X, \mathcal M)$ generated by
the polynomials $x\la ijk \ra$ ($1\leq i<j<k\leq 6$).
\end{lemma}
Let $\widehat{pl}$ be the map from $\Mp(3,6)$ to $\C^{10}$ defined by 
$$
\widehat{pl}: \Mp(3,6)\ni x \mapsto 
(\dots,x\la J\ra,\dots)_{\la J\ra \in P_{3,3}}\in \C^{10},
$$
where we arrange $x\la J \ra$ lexicographically
for $J=\{j_1,j_2,j_3\}$ with $j_3\le 5$. 
By Lemma \ref{Pl-rel}, the image of $pl$ is contained in a $5$-dimensional linear 
subspace of $\C^{10}$. The map $X\to \P^4$ induced from $\widehat{pl}$ is
denoted by $pl$.

The space of relative global holomorphic $2$-forms
$H^0(\widehat{\cS}, \Omega_{\widehat{\cS}/\Mp(3,6)}^2)$ 
is generated by
$$
\varphi={t_0dt_1\wedge dt_2-t_1dt_0\wedge dt_2+t_2dt_0\wedge dt_1
\over y}.
$$\begin{proposition}
\label{invariant}
The form $\varphi$ satisfies the equality:
$$
(g,\lambda_i,\lambda)^*\varphi=\rho^{-1}(g,\lambda_i,\lambda)\varphi.
$$
Therefore it defines a global section of 
$H^0(\cS, \Omega_{\cS/X}^2\otimes \mathcal L^{-1})$.

\end{proposition}

\subsection{Topological cycles at a reference point}
\label{topological cycle at ref pt}
We take a reference point $\dot x$
in $\Mp(3,6)$ as 
$$\dot x=\pmatrix{
p_1^2& p_2^2& p_3^2& p_4^2& p_5^2& p_6^2\cr
-p_1& -p_2& -p_3& -p_4& -p_5& -p_6\cr
1 & 1 & 1 &1 & 1 & 1\cr
}$$
with $p_1=-3,\ p_2=-2,\ p_3=-1,\ p_4=1,\ p_5=2,\ p_6=3.$ We consider
the affine
coordinates $s_1=t_1/t_0$ and $s_2=t_2/t_0$ of $\P^2$.
We construct topological 2-cycles of $\cS(\dot x)$
using the isomorphism of the Kummer surface of $C$
and $\cS(\dot x)$ given in \cite{T}.

Let $C$ be a hyperelliptic curve defined by
$$
u^2=\prod_{i=1}^6(w-p_i),
$$
and $C_1, C_2$ be copies of $C$.
Let $sym:C_1\times C_2\to \cS(\dot x)$ be a map defined by
$$
((w_1,u_1),(w_2,u_2))\mapsto (s_1,s_2,y)=(w_1+w_2,w_1w_2,u_1u_2).
$$
For $a<b \in \R$, the 1-chain in $\P^1$
defined by the segment from $a$ to $b$ is
denoted by $(a,b)$. We define chains $A_1', A_2', B_1'$ and $B_2'$ in $C$
by the liftings of $(p_1,p_2), (p_5,p_6), (p_2,p_3)$ and $(p_4,p_5)$
on which $1/u$ is in $\i\R_+,\i\R_+,\R_+$ and $\R_+$.
Then $A_i=A_i'-\sigma(A_i')$ and $B_i=B_i'-\sigma(B_i')$
become 1-cycles on $C$.

We set $A_1=\gamma_1,A_2=\gamma_2,B_1=\gamma_3$ and $B_2=\gamma_4$.
We define a topological cycle ${\overline{\gamma}_{ij}}^{*}$ by 
$sym_*(\overline{\gamma}_{i}\times
\overline{\gamma}_{j}).$ The proper inverse image of ${\overline{\gamma}_{ij}}^{*}$ in
$\cS(\dot x)$ is denoted by $\overline{\gamma}_{ij}$. Then $\gamma_{ij}=\overline{\gamma}_{ij}/2$
is an element in
$H_2(\cS(\dot x),\Z)$. 
Let $\g'_{ij}$ be the orthogonal projection of $2\g_{ij}$ to $T(x)$.
Then
$\{\g'_{12},\g'_{13},\g'_{14},\g'_{23},\g'_{24},\g'_{34}\}$
is a basis of
$T(x)$.
Since
$$
sym^*(\f(\dot x))=
{(w_1-w_2)dw_1\wedge dw_2\over 
u_1u_2
},
$$
we have
\begin{equation}
\label{arg}
\displaystyle
\iint_{\g_{12}}\f(\dot x)\in-\Rp,
\quad \iint_{\g_{34}}\f(\dot x)\in \Rp, 
\end{equation}
$$
\iint_{\g_{13}}\f(\dot x),\ 
\iint_{\g_{14}}\f(\dot x)\in \i\Rp,\quad 
\iint_{\g_{23}}\f(\dot x),\ \iint_{\g_{24}}\f(\dot x)\in -\i\Rp.
$$

\begin{proposition}
\label{cycles}
We set
\begin{eqnarray*}
\g=\tr(\g_{12},\g_{13},\g_{14},\g_{23},\g_{24},\g_{34}),\quad
\g'=\tr(\g'_{12},\g'_{13},\g'_{14},\g'_{23},\g'_{24},\g'_{34}).
\end{eqnarray*}
Then the intersection matrix is equal to
$$(\g \cdot \tr\g')=H,\quad (\g' \cdot \tr\g')=2H,$$
where 
$$
H=-\pmatrix{
 & & & & & 1\cr 
 &1& & & &  \cr 
 & & &1& & \cr 
 & &1& & & \cr 
 & & & &1& \cr 
1& & & & & \cr }.
$$
Thus the lattice structure of $T(x)$ is
equal to $U(2)\oplus U(2) \oplus A_1(-1) \oplus A_1(-1)$.
\end{proposition}
(For details, 
see \cite{MSY} and \cite{Y}, Chapter VIII.)

For any $x\in \Mp(3,6)$, take a path $\rho_x$ in $\Mp(3,6)$ 
connecting $\dot x$ and $x$, and define bases $\g_{ij}(x)$ and $\g'_{ij}(x)$ 
as the continuations of $\g_{ij}(\dot x)$ and 
$\g'_{ij}(\dot x)$ along the path $\rho_x$ by the local triviality.
They depend only on the homotopy class of $\rho_x$.
The intersection matrix for 
$(\g'(x) \cdot \tr\g'(x))$
are equal to that in Proposition \ref{cycles}.
Let $\w$ be a vector defined by
\begin{equation}
\label{periods}
\w(x)=\tr(\w_{12}(x),\w_{13}(x),\w_{14}(x),\w_{23}(x),\w_{24}(x),\w_{34}(x)),
\end{equation}
$$
\w_{ij}(x)=\iint _{\g_{ij}(x)}\f(x).
$$
The map
\begin{equation}
\label{eq:period map first}
\widetilde{per}:\widetilde{X}\ni x\mapsto [\w(x)]\in\P^5,
\end{equation}
is called the period map, where $\widetilde{X}$ is the universal
covering of $X$.

\subsection{Period integrals and Hypergeometric functions}


We define two hypergeometric series $F_S^\a(z)$ and $F_T^\a(z)$ 
of variables $z=\pmatrix{z_1& z_3 \cr z_2 &z_4}$ 
with parameters $\a=(\a_1,\dots,\a_6)$ satisfying $\sum_{j=1}^6\a_j=3$
as 
\begin{eqnarray}
\label{S-series}
F_S^\a(z)&=&\sum_{n\in \N^4}
{(1\!-\!\a_1)_{n_1+n_3}(1\!-\!\a_2)_{n_2+n_4}(\a_5)_{n_1+n_2}(\a_6)_{n_3+n_4}
\over (2\!-\!\a_1\!-\!\a_3)_{n_1+n_3}
(2\!-\!\a_2\!-\!\a_4)_{n_2+n_4}n_1!n_2!n_3!n_4!}
z^n,\\
\label{T-series}
F_T^\a(z)&=&\sum_{n\in \N^4}
{(1-\a_1)_{n_1+n_3}(1-\a_2)_{n_2+n_4}(\a_5)_{n_1+n_2}(\a_6)_{n_3+n_4}
\over (3-\a_1-\a_2-\a_3)_{n_1+n_2+n_3+n_4}n_1!n_2!n_3!n_4!}
z^n,
\end{eqnarray}
where $\N=\{0,1,2,\dots\}$, 
$z^n=z_1^{n_1}z_2^{n_2}z_3^{n_3}z_4^{n_4}$ for $n=(n_1,\dots,n_4)$, 
 $(\a_j)_{n_j}=\a_j(\a_j+1)\cdots(\a_j+n_j-1)
=\G(\a_j+n_j)/\G(\a_j)$, 
and we assume that 
$$
\begin{array}{cll}
\a_1+\a_3-2,\ \a_2+\a_4-2\notin \N &\textrm{for} &F_S^\a(z),\\
\a_1+\a_2+\a_3-3\notin \N &\textrm{for} &F_T^\a(z).
\end{array}
$$ 
They absolutely 
converge on the domain 
$$\{z\in \C^4\mid |z_1|+|z_2|<1,\ |z_3|+|z_4|<1\}.$$
By the standard argument for Euler type integrals, we have the following Proposition.
\begin{proposition}
\label{int-rep}
The hypergeometric series $F_S^\a$ and $F_T^\a$ 
admit the integral representations:
\begin{eqnarray*}
F_S^\a(z)&=&{1\over B(1-\a_1,1-\a_3)B(1-\a_2,1-\a_4)}
\int_0^1\!\!\!\int_0^1 L_S^\a(z,s)\ ds_1ds_2,\\
F_T^\a(z)&=&
{\G(3-\a_1-\a_2-\a_3)\over \G(1-\a_1)\G(1-\a_2)\G(1-\a_3)}
\iint_\De L_T^\a(z,s) ds_1ds_2,
\end{eqnarray*}
where $B$ denotes the beta function,
\begin{eqnarray*}
L_S^\a(z,s)&=&
s_1^{-\a_1}s_2^{-\a_2}(1-s_1)^{-\a_3}(1-s_2)^{-\a_4}\\
& &\times (1-z_1s_1-z_2s_2)^{-\a_5}(1-z_3s_1-z_4s_2)^{-\a_6},\\
L_T^\a(z,s)&=&
s_1^{-\a_1}s_2^{-\a_2}(1-s_1-s_2)^{-\a_3}\\ 
& &\times (1-z_1s_1-z_2s_2)^{-\a_5}(1-z_3s_1-z_4s_2)^{-\a_6},\\
\De&=&\{(s_1,s_2)\in \R^2\mid s_1>0,s_2>0,s_1+s_2<1\},
\end{eqnarray*}
$\arg(s_j)=\arg(1-s_j)=\arg(1-s_1-s_2)=0$ on 
each interior area of the integrations, 
$\arg(1-z_1s_1-z_2s_2)$, $\arg(1-z_3s_1-z_4s_2)$ become $0$ 
at $(s_1,s_2)=(0,0)$, 
and we assume that 
$$\begin{array}{cll}
\re(\a_1),\re(\a_2),\re(\a_3),\re(\a_4)<1 &\textrm{for} &F_S^\a(z),\\
\re(\a_1),\re(\a_2),\re(\a_3)<1 &\textrm{for} &F_T^\a(z).
  \end{array}
$$
\end{proposition}

From now on, we put $\a=(1/2,\dots,1/2)$ and we set
$$F_S(z)=F_S^\a(z),\quad F_T(z)=F_T^\a(z).$$
We can regard $ x\la J\ra\w_{ij}( x)^2$ 
as a multivalued function on $X$ for any $\la J\ra$.
Proposition \ref{invariant} and \ref{int-rep}
imply the following.
\begin{proposition} 
\label{series-exp}

For a fixed  $\la J\ra\in P_{3,3}$,
the product $ x\la J\ra\cdot\w_{ij}( x)^2$ is 
invariant under the action of 
$GL_3(\C)\times \T$, i.e., 
$$
(g, \lambda_i,\lambda)^*(x\la J\ra \cdot \w_{ij}^2)
=x \la J\ra \cdot \w_{ij}^2.
$$
As a consequence,  
if $x$ is in a neighborhood of $\dot x\in \Mp(3,6)$, 
$\w_{ij}(x)^2$ can be expressed as 
\begin{eqnarray*}
\w_{ij}(x)^2x\la J\ra=
\left\{
\begin{array}{ll}
 4\pi^4F_S(\zeta_{ij}(x))^2{\nu_{ij}\la J\ra  }, &
 (i,j)=(1,2), (1,4), (2,3), (3,4),\\[2mm]
 16\pi^2F_T(\zeta_{ij}(x))^2{\nu_{ij}\la J\ra }, &
  (i,j)=(1,3), (2,4),
\end{array}
\right.
\end{eqnarray*}
where we set $\zeta_{ij}(x)=\pmatrix{z_1& z_3\cr z_2& z_4\cr}$ 
so that the following $\nu_{12},\dots,\nu_{34}\in \Mp(3,6)$ 
are equivalent to $x$ as elements of $X$: 
{\small
$$\matrix{
\nu_{12}=\pmatrix{
 0 & 1 & 1 & 1 & 1 & 0\cr
 1 &-1 &-z_1 &-z_3 & 0 &0\cr
 0 & 0 &-z_2 &-z_4 & -1&1 \cr},& 
\nu_{13}=\pmatrix{
0 & 0 & 1 & 1 & 1 & 1 \cr
0 & 1 &-1 &-z_1 &-z_3 & 0 \cr
1 & 0 &-1 &-z_2 &-z_4 & 0 \cr}, 
\cr
\nu_{14}=\pmatrix{
 1 & 0 & 1 & 0 & 1 & 1\cr
 -1 &1 &-z_1 & 0 & 0 &-z_3\cr
 0 & 0 &-z_2 & 1 & -1&-z_4 \cr},&
\nu_{23}=\pmatrix{
 1 & 1 & 0 & 1 & 1 & 0 \cr
-z_1 &0& 0 &-z_3& -1 & 1\cr
-z_2&-1& 1 &-z_4& 0 &0 \cr}, 
\cr
\nu_{24}=
\pmatrix{
1 & 1 & 1 & 1 & 0 & 0 \cr
0 &-z_1 &-z_3 & -1& 1 &0 \cr
0 &-z_2 &-z_4 & -1& 0 &1 \cr}, &
\nu_{34}=\pmatrix{
 1 & 1 & 0 & 0 & 1 & 1 \cr
-z_1 &-1& 1 & 0& 0 &-z_3 \cr
-z_2 &0& 0 & 1& -1 &-z_4 \cr}. 
}
$$
}

\end{proposition}

\subsection{Preparation for the association involution}
\label{association preparation}
Let $as$ be an automorphism of $\Mp(3,6)$ given by 
$$as:\Mp(3,6)\ni (y_1,y_2)\mapsto (\tr(y_1^{-1} y_2y_1),\tr y_1)\in 
\Mp(3,6),$$ 
where $y_1,y_2 \in GL_3(\C)$.
By a straightforward calculation, we have 
$$
\widetilde{pl}\circ as=\widetilde{pl},\quad 
as^2(y_1,y_2)=(y_1^{-1}y_2y_1y_2^{-1}y_1,y_1^{-1}y_2y_1)
=y_1^{-1}y_2y_1y_2^{-1}(y_1,y_2).
$$ 
Therefore $as$ induces an involution on $X$, which is called the
association involution and also denoted by $as$.
By the above equality, we have an induced morphism
$$
pl^*:X/\la as\ra \to \P^4.
$$

\begin{proposition}[Chapter VII of \cite{Y}] 
The morphism $pl^*$
is an open immersion.
\end{proposition}
Let $\overline{X}$ be the normalization of $\P^4$ in $X$.
Then we have the diagram:
$$
\begin{array}{ccc}
X & \to & \overline{X} \\
\downarrow & & \downarrow \\
X/\la as \ra &\to & \P^4.
\end{array}
$$
The induced map $\overline{X}\to \P^4$ is denoted as $\overline{pl}$.
Let $x=(x_{ij})_{ij} \in \Mp(3,6)$. 
We define the following
polynomials
\begin{eqnarray*}
Q&=&\det
 (x_{1i}^2,x_{2i}^2,x_{3i}^2,x_{2i}x_{3i},
x_{3i}x_{1i},x_{1i}x_{2i})_{i=1,\dots, 6}, \\
D(ijk)&=&\det(x_{pi},x_{pj},x_{pk})_{p=1,\dots, 3}, \\
\{ij;kl\}&=&D(ijm)D(ijn)D(mkl)D(nkl), \\
T(ijklmn)&=&D(ijk)D(klm)D(mni)D(nlj)
\end{eqnarray*}
 for 
$\{i,j,k,l,m,n\}=\{1,\dots,6\}$.
Then we have $\{ij;kl\}=\pm x\la ijm \ra x\la ijn\ra$ and 
$as(Q)=-Q$.
We give an explicit description of the normal form
$$
\left(
\begin{array}{cccccc}
1 & 0 & 0 & 1 & 1 & 1  \\
0 & 1 & 0 & 1 & x_1 & y_1  \\
0 & 0 & 1 & 1 & x_2 & y_2  
\end{array}
\right)
$$
of the inverse image under $\overline{pl}$.
By Lemma A6.8 of \cite{MSY},
we have
$$
x_2=\frac{D(125)D(234)}{D(235)D(124)}=
\frac{T(125364)}{\{35;14\}}.
$$
We have $as(T(125364))=T(364125)$ and
by Lemma A7.3 of \cite{MSY},
we have
\begin{eqnarray*}
T(125364)+as(T(125364))&=&
\{14;53\}-\{52;16\}+\{63;54\} \\
& &-\{23;15\}+\{24;56\}\\
T(125364)\cdot as(T(125364))
&=&\{16;23\}\cdot \{12;36\}
\end{eqnarray*}
Therefore $T(125364)$ is defined by a quadratic equation with the coefficients
in polynomials of $x\la ijk \ra$'s.
The values $x_1, y_2, y_2$ can be obtained by substituting
indices $2 \leftrightarrow 3$, $5\leftrightarrow 6$.

Since
$$
T(125364)-as(T(125364))=Q,
$$
if $\{16;23\}\cdot \{12;36\}=0$,
then $Q=T(125364)$ or $-T(364125)$ and the values
$x_2, x_1, y_2, y_1$ are polynomials of $x\la ijk \ra$.
Thus we have the following lemma.
\begin{lemma}
The inverse image of the divisor 
$\{x\la 164 \ra=0\}$ on $\P^4$
under the map $\overline{pl}$
consists of two irreducible components.
On $\overline{X}$, we can express 
$x_1,x_2,y_1, y_2$ as rational functions of
$x\la ijk\ra$ on each irreducible component.
\end{lemma}
For an explicit description of the inverse image, see Proposition \ref{x123=0}.

\section{Bounded symmetric domains and theta functions.}

\subsection{Period map and symmetric domains $D_H$ and $\D$}
\label{isom symmetric domain}
In this section, we introduce two symmtric domains $D_H$ and $\D$.
The target space for the natural period map for \KK surfaces 
is the symmetric domain $D_H$ of type $IV$. We use theta functions on
$\D$ and an isomorphism between $D_H$ and $\D$ to construct automorphic functions on 
$D_H$ by using results in \cite{Ma}.

By Riemann bilinear relations for the \KK surface $\cS(x)$, 
and the choice of orientations of $T(x)$ in (\ref{arg}),
the class 
$[\w(x)]$ of (\ref{periods}) in $\P^5$
belongs to the subset
$$
D_H =\{[w]\in \P^5   \mid \tr\w\; H\;  \w=0,\quad
\w^*\;H \;\w>0, \quad
\im({\w_{14}\over \w_{34}})
>0
\}.
$$
Here $y^*=\tr \overline y$ denotes the adjoint of a matrix $y$.
Therefore we regard the map (\ref{eq:period map first})
as $\widetilde{per}:\widetilde{X}\to D_H$.

Let $\D$ be the bounded symmetric domain of type $I_{22}$ defined by
$$
\D=\{\tau \in M(2,2)\mid {\tau-\tau^*\over 2\i}\textrm{ is positive
definite } \}.
$$ 
In this subsection, 
we define an isomorphism $\D \to D_H$.
Let $\tau$ be an element in $\D$. We set 
$$\tilde \tau=\pmatrix{\tau\cr E_2}.$$
Let $\tilde\tau\la i_1i_2\ra$ be the
$(i_1,i_2)\times (1,2)$-minor of the $4\times 2$ matrix $\tilde\tau$. 
They satisfy the Pl\"ucker relation
$$
\tilde\tau\la12\ra\tilde\tau\la34\ra
-\tilde\tau\la13\ra\tilde\tau\la24\ra
+\tilde\tau\la14\ra\tilde\tau\la23\ra=\tr v(\tilde\tau)\ H'\ v(\tilde\tau)
=0,
$$
where 
$$H'=\pmatrix{
 & & & & & 1\cr 
 & & & &-1&  \cr 
 & & &1& & \cr 
 & &1& & & \cr 
 &-1& & & & \cr 
1& & & & & \cr },\quad 
v(\tilde\tau)=
\pmatrix{\tilde\tau\la12\ra\cr \tilde\tau\la13\ra\cr 
\tilde\tau\la14\ra
\cr \tilde\tau\la23\ra\cr \tilde\tau\la24\ra\cr 
\tilde\tau\la34\ra}.$$
Since the matrix ${(\tau-\tau^*)/ 2\i}$ is positive definite, we have
$$v(\tilde\tau)^*\; H\; v(\tilde\tau)>0,\quad \im({\tilde\tau\la14\ra(\tau)\over \tilde\tau\la34\ra(\tau)})>0.$$
We set $$Q=\pmatrix{1 & & & & & \cr
&{1+\i\over 2}&&& {-1+\i\over 2}&\cr
&&1&&&\cr
&&&1&&\cr
&{-1+\i\over 2}&&&{1+\i \over 2}&\cr
&&&&&1\cr
}.$$
Then by the equality
$$
Q^*\; H\; Q=H, \quad H'=\tr Q\; H\; Q,\quad H=\tr Q\; H'\; Q,
$$
the class $[Qv(\tilde\tau)]$ of $Qv(\tilde\tau)$
in $\P^5$ is contained in $D_H$. 
Thus we have an isomorphism
\begin{equation}
\label{fundmental isom on domains}
p\jmath_\D:\D \ni\tau \mapsto  [Qv(\tilde\tau)] \in D_H.
\end{equation}
We define the normalized period 
matrix of $\cS(x)$ 
by $\tau=\tau(x)=p\jmath_\D^{-1}([\omega(x)])\in \D$.
Then we have
\begin{equation}
\label{periodmatrix}
\tau(x)={1\over \w_{34}(x)}
\pmatrix{\w_{14}(x)& -{\w_{13}(x)-\i\w_{24}(x)\over1+\i} \cr 
-{\w_{13}(x)+\i\w_{24}(x)\over1-\i}& -\w_{23}(x)},
\end{equation}
where $\w_{ij}(x)$ are defined by (\ref{periods}).
Thus we have a map 
$$
\widetilde{per}:\widetilde{X}\ni x\mapsto \tau(x)\in \D.$$
Since $\tilde\tau\la34\ra(\tau)=1$, we have 
\begin{equation}
\label{ItoIV}
\jmath_\D(\tau(x))=\w(x)/\w_{34}(x)
\end{equation}
as elements of $\C^6$.

\subsection{Homomorphisms of discrete groups}
We define
a discrete group $\bG_H$ in $GL_6(\Z)$ by
$$\bG_H=\{R\in GL_6(\Z)\mid \tr R H R=H,\ 
\im({[R\w(\dot x)]_{14}\over[R\w(\dot x)]_{34}})>0\},$$
where $[R\w(\dot x)]_{ij}$ denotes the $(ij)$-component of $R\w(\dot x)$.
Its center consists of $\pm E_6$. The the group $\bG_H$ acts on $D_H$
from the left.
Since the monodromy action preserves the intersection forms,
the monodromy group is contained in $\bG_H$.

We define $U_{22}(\Z[\i])$ and the principal congruence subgroup of 
level $(1+\i)$ by 
$$U_{22}(\Z[\i])=\left\{g\in GL_4(\Z[\i])\; \Big|\; 
g\; I_{22}\; g^*=I_{22}
\right\},$$
$$
U_{22}(1+\i)=\{g\in U_{22}(\Z[\i])\mid g\equiv E_4 \bmod (1+\i)\}.
$$
Then an element 
$g=\pmatrix{g_{\underline{11}} & g_{\underline{12}}\cr
g_{\underline{21}} & g_{\underline{22}}\cr}$
in
$U_{22}(\Z[\i])$ 
acts on $\D$ by
$$g\cdot \tau=(g_{\underline{11}} \tau + g_{\underline{12}})
(g_{\underline{21}}\tau + g_{\underline{22}})^{-1},$$
where 
$I_{22}=\pmatrix{O & -E_2\cr E_2 &O}$ and 
$g_{\underline{ij}}\in M(2,2)$.

In this subsection, we define a homomorphism 
$U_{22}(\Z[\i])\to \bG_H/\la \pm 1\ra$ 
of discrete group which is compatible with the isomorphism
of symmetric domains $\D \to D_H$ defined in 
(\ref{fundmental isom on domains}).
For $g=(g_{ij})=\pmatrix{g_{\underline{11}} & g_{\underline{12}}\cr
g_{\underline{21}} & g_{\underline{22}}}\in U_{22}(\Z[\i])$, 
we set a $6\times 6$ matrix $\wedge ^2 g$ by  
$$\wedge ^2 g=\left(\left|\matrix{
g_{i_1j_1} &g_{i_1j_2}\cr
g_{i_2j_1} &g_{i_2j_2}\cr}\right|
\right)_{(i_1i_2),(j_1j_2)},
$$
where $1\le i_1<i_2\le 4$, $1\le j_1<j_2\le 4$, and  
they are arranged lexicographically. Then we have
$$\det(\wedge ^2 g)=\det(g)^3,$$
and 
$$\det(g_{\underline{21}}\tau+g_{\underline{22}})
v(\tilde \tau')=v(g\tilde\tau) =(\wedge^2g)v(\tilde\tau)$$
as elements of $\C^6$, where $\tau\in \D$,
$$\tau'=g\cdot \tau=(g_{\underline{11}} \tau + g_{\underline{12}})
(g_{\underline{21}}\tau + g_{\underline{22}})^{-1},
\quad \tilde\tau'=\pmatrix{\tau'\cr E_2\cr},\quad 
g\tilde\tau=g\pmatrix{\tau\cr E_2\cr}.$$
Thus we have 
$$
\det(g_{\underline{21}}\tau+g_{\underline{22}})
Q v(\tilde\tau')=\{Q(\wedge^2g)Q^{-1}\} \{Q v(\tilde\tau)\},$$
and by the definition of $\jmath_\D$, we have
\begin{equation}
\label{exptau}
\det(g_{\underline{21}}\tau+g_{\underline{22}})\jmath_\D(g\cdot\tau)
=Q(\wedge^2g)Q^{-1}\jmath_\D(\tau).
\end{equation}
The matrix $Q(\wedge^2g)Q^{-1}$ belongs to the orthogonal group 
with respect to the quadratic form $H$.
Moreover, a straight forward calculation shows 
$$Q(\wedge^2g)Q^{-1}\in \left\{
\begin{array}{rll}
 SL_6(\Z) & \textrm{if} & \det(g)=1,\\
\i\ SL_6(\Z) & \textrm{if} & \det(g)=-1,
\end{array}
\right.
$$
for $g\in U_{22}(\Z[\i])$. 
We set
\begin{equation}
\label{4to6}
R_{g}=\sqrt{\det(g)} Q(\wedge^2g)Q^{-1},
\end{equation}
which is determined modulo sign. 
Since  $\wedge^2(\i g)=-\wedge^2 g$, 
we have $R_{\i g}=-R_g$. 
Thus $R_g$ defines a homomorphism 
$$
\begin{array}{ccc}
U_{2,2}(\Z[\i])/\la iE_4 \ra  & \to &  
\bG_H/\la \pm E_6\ra \\
\cap & & \cap \\
Aut(\D) & \to & Aut(D_H) 
\end{array}.
$$

\subsection{Monodromy actions on the spaces $D_H$ and $\D$}

Each center of $U_{22}(Z[\i])$ and $U_{22}(1+\i)$ is 
the group $\la \i E_4\ra$ generated by the scalar matrix $\i E_4$.
We have 
$$\tr(\bar g \cdot \tau)=g\cdot \tr\tau$$
for any $g\in U_{22}(\Z[\i])$ and $\tau\in \D$.
Let $tp$ be the transpose operator acting on $\D$ and 
$\la tp\ra$ be the group generated by $tp$. 
The fixed locus of $tp$ is the Siegel upper half space 
$\H_2=\{\tau\in \D\mid \tr\tau=\tau\}$ of degree $2.$
We define $U_{22}^{tp}(\Z[\i])$ acting on $\D$  as  the group 
generated by $U_{22}(\Z[\i])/\la \i E_4\ra$  and $\la tp\ra $ with relations 
$$ (tp) g= \bar g (tp)$$
for any $g\in  U_{22}(\Z[\i])$.
This group is a semi-direct product  
$(U_{22}(\Z[\i])/\la \i E_4\ra) \rtimes 
\la tp\ra$.
We set $U_{22}^{tp}(1+\i)=(U_{22}(1+\i)/\la \i E_4\ra)\rtimes \la tp\ra$.


\begin{proposition}[\cite{Ma},\cite{KiM},\cite{Y}]
\label{monodromy}
\begin{itemize}
\item[$(1)$]

We define
the principally congruence subgroup $\bG_H(2)$ 
of level $2$ by
\begin{eqnarray*}
\bG_H(2)=\{R\in \bG_H \mid  R\equiv E_6\bmod 2
\}.
\end{eqnarray*}
Then the monodromy group for $\widetilde{per}:\widetilde{X} \to D_H$
is equal to $\bG_H(2)$.
\item[$(2)$]
The monodromy group of $\widetilde{per}:\widetilde{X}\to \D$ over
$X$ is equal to
$$U_{22}^M(1+\i)=\{(g,tp^k)\in U_{22}^{tp}(1+\i)
\mid \det(g)=(-1)^k\}.$$
We note that $\det (g)$ is well defined on $U_{22}(\Z[\i])/\la \i E_4\ra$.
The monodromy group over $X/\la as\ra$
is equal to $U_{22}^{tp}(1+\i)$.
The map $per\circ {pl^*}^{-1}$ induces the isomorphism from 
$\P^4$ to the Satake compactification of the quotient $\D/U_{22}^{tp}(1+\i)$. 
\item[$(3)$]
Let $As$ be an element in $GL_6(\Z)$ defined by
$$
As(\w_{ij})=\left\{
\begin{array}{ll}
\w_{ij} &(i,j) \neq (1,3), (2,4) \\
\w_{13} &(i,j) = (2,4)\\
\w_{24} &(i,j) = (1,3),
\end{array}\right.
$$
and $as$ the association involution defined in 
\S \ref{association preparation}. 
Then we have
\begin{equation}
\label{assoc and T}
\w(as(x))=As(\w(x))
\end{equation}
for $x$ in a small neighborhood of $\dot x\in X$. 

\item[$(4)$]
Under the isomorphism $Aut(D_H)\simeq Aut(\D)$, the matrix $As$ defined in (3) corresponds to $tp$.
Therefore we have
\begin{equation}
\label{as-tr}
\tau(as(x))=tp(\tau(x))
\end{equation}
for $x$ in a small neighborhood of $\dot x\in X$ and
an isomorphism
\begin{equation}
\label{isom of discrete group extended}
\widetilde \bG_H(2)\simeq U_{22}^{tp}(1+\i),
\end{equation}
where $\widetilde\bG_H(2)=\bG_H(2)/\la \pm E_6\ra\cdot \la As\ra$.

\end{itemize}
\end{proposition}
By the above proposition, we have a map
$$
per:X\to D_H/\bG_H(2) \simeq \D/U^{M}_{22}(1+\i).
$$
Since the last components $\jmath_\D(\tau)$ and $\jmath_\D(g\cdot \tau)$
are $1$, we have
$$
\pm \det(g_{\underline{21}}\tau(x) + g_{\underline{22}})
=\pm {\sqrt{\det(g)}[R_g\w(x)]_{34}\over \w_{34}(x)}
$$
by the equality (\ref{exptau}) together with 
(\ref{ItoIV}) and (\ref{4to6}),
where  $[R_g\w(x)]_{34}$ denotes the 
$(34)$-component of the column vector $R_g\w(x)$.
By squaring this equality, we have 
\begin{equation}
\label{auto-factor}
\det (g)\det(g_{\underline{21}}\tau(x) + g_{\underline{22}})^2
={[R_{g}\w(x)]_{34}^2\over \w_{34}(x)^2}
\end{equation}
for $g\in U_{22}(\Z[\i])/\la \i E_4\ra$.

\subsection{Theta functions and their functional equations}

The theta function $\hh_{ab}$ with characteristic $a,b$ on $\D$ is defined as
\begin{equation}
\label{Theta on D}
\hh_{ab}(\tau)=\sum_{n\in \Z[\i]^2}
\ex[\frac{1}{2}(n+a)\tau(n+a)^*+\re((n+a)b^*)],
\end{equation}
where $x^*=\tr \bar x$, $\ex[x]=\exp(2\pi \i x)$, $\tau\in \D$, 
$n=(n_1,n_2)\in \Z[\i]^2$, $a=(a_1,a_2),b=(b_1,b_2)\in \Q[\i]^2$.

\begin{remark}
\label{decompose}
This $\hh_{ab}$ is different from 
that defined in \cite{MY} and \cite{Ma} by the factor of $\ex[\re(ab^*)]$ .
If $\tau$ belongs to the Siegel upper half space $\H_2$ of degree $2$, then
$\hh_{ab}$ decomposes into the product of Riemann's theta constants:
$$\hh_{ab}(\tau)=\vartheta_{\re(a)\re(b)}(\tau)\vartheta_{\im(a)\im(b)}(\tau),
$$
where 
$$\vartheta_{a'b'}(\tau)=\sum_{n\in \Z^2}
\ex[\frac{1}{2}(n+a')\tau\tr(n+a')+(n+a')\tr b']
$$
for $a',b'\in \Q^2$.
\end{remark}

This function satisfies 
$$\hh_{ab}(\tr\tau)=\hh_{\bar a,\bar b}(\tau),\quad 
\hh_{a+n,b}(\tau)=\hh_{ab}(\tau),\quad 
\hh_{a,b+n}(\tau)=\ex[\re(an^*)]\hh_{ab}(\tau),$$
for any $n\in \Z[\i]^2$.
For $a,b\in \Z[\i]^2$,  
$\hh_{{a\over 1+\i}{b\over 1+\i}}(\tau)$ is denoted by 
$\hh_{[ab]}(\tau)$. 
Then
$\hh_{[ab]}^2(\tau)$ depends only on the class  
of $a$ and $b$ in $\F_2^2\simeq (\Z[\i]/(1+\i)\Z[\i])^2$.
We set
$$
Ev=\{(a,b) \in (\Z[\i]/(1+\i)\Z[\i])^2 \mid ab^*=0 \textrm{ mod }(1+\i)\}.
$$
Then we have $\# Ev=10$ and
$\hh_{[ab]}(\tau)=0$ if $(a,b)\notin Ev$.  
We identify the sets $P_{3,3}$ and $Ev$ by the rule given in 
Table \ref{part and char}.
Under this correspondence,  
$\hh_{[ab]}^2(\tau)$ is denoted by $\hh_{\la J\ra}^2(\tau)$. 

\begin{table}
$$\begin{array}{c}
\la 123\ra \leftrightarrow [1111]\\
\la 124\ra \leftrightarrow [0011]\\
\la 125\ra \leftrightarrow [0010]\\
\la 134\ra \leftrightarrow [0001]\\
\la 135\ra \leftrightarrow [0000]\\
 \end{array}
\qquad
\begin{array}{c}
\la 145\ra \leftrightarrow [1100]\\
\la 234\ra \leftrightarrow [1001]\\ 
\la 235\ra \leftrightarrow [1000]\\
\la 245\ra \leftrightarrow [0100]\\
\la 345\ra \leftrightarrow [0110]\\
\end{array}
$$
\caption{Correspondence between $\la J\ra $ and $[ab]$}
\label{part and char}
\end{table}
\begin{remark}
The correspondence between $P_{3,3}$ and $Ev$ is different from that in 
\cite{Ma}, since the bases of the transcendental lattice $T(x)$ are different.
\end{remark}
\begin{proposition}[\cite{Ma}]
\label{hirei}
\begin{itemize}
\item[$(1)$]
They satisfy 
\begin{eqnarray*}
\hh_{[ab]}^2(\tr\tau)&=&\hh_{[ab]}^2(\tau),\\
\hh_{[ab]}^2(g\cdot \tau)&=&
\det(g)\det(g_{\ul{21}}\tau +g_{\ul{22}})^2\hh_{[ab]}^2(\tau),
\end{eqnarray*}
for any $g=\pmatrix{g_{\ul{11}}&g_{\ul{12}}\cr g_{\ul{21}}&g_{\ul{22}}}
\in U_{22}(1+\i)$.

\item[$(2)$]
We define a map $\theta:\D/U_{22}^{tp}(1+\i) \to \P^9$ by
$$
\tau
\mapsto[\dots,\hh_{\la J\ra}^2(\tau),\dots]_{\la J\ra\in P_{3,3}}\in\P^9.
$$ 
Then the map $pl$ is equal to the composite $\theta\circ per$ from $X$ to
$\P^9$.
\end{itemize}
\end{proposition}
\begin{theorem}[$2\tau$-formula]
\label{2tau-formula}
We have 
$$4\hh_{ab}(2\tau)=\sum_{q\in \F_2^2}
\ex[-\re(aq^*)]\hh_{(1+\i)a,{b+q\over 1-\i}}(\tau),$$
where $q$ runs over the set $\F_2^2=\{(0,0), (0,1), (1,0), (1,1)\}$. 
In particular, 
\begin{eqnarray*}
\hh_{[0000]}(2\tau)&=&{1\over 4}
(\hh_{[0000]}(\tau)+\hh_{[0001]}(\tau)+\hh_{[0010]}(\tau)+\hh_{[0011]}(\tau)),\\
\hh_{[0100]}(2\tau)&=&{1\over 4}
(\hh_{[0000]}(\tau)-\hh_{[0001]}(\tau)+\hh_{[0010]}(\tau)-\hh_{[0011]}(\tau)),\\
\hh_{[1000]}(2\tau)&=&{1\over 4}
(\hh_{[0000]}(\tau)+\hh_{[0001]}(\tau)-\hh_{[0010]}(\tau)-\hh_{[0011]}(\tau)),\\
\hh_{[1100]}(2\tau)&=&{1\over 4}
(\hh_{[0000]}(\tau)-\hh_{[0001]}(\tau)-\hh_{[0010]}(\tau)+\hh_{[0011]}(\tau)).
\end{eqnarray*}
\end{theorem}
\pr
We consider the summation 
\begin{equation}\label{sum-L}
\sum_{n'\in L}\sum_{q\in \F_2^2}
\ex[\frac{1}{2}\re(n'q^*)]\cdot\ex[(n'+a)\tau(n'+a)^*+\re((n'+a)b^*)],
\end{equation}
where $L={1\over 1+\i}\Z[\i]^2$. 
Since $\Z[\i]^2\subset L$, $L/\Z[\i]^2\simeq \F_2^2$ and 
$\ex[\frac{1}{2}\re(n'q^*)]$ $(q\in \F_2^2)$ are the characters of 
the quotient group $L/\Z[\i]^2$,  
this summation  reduces to the four times of 
the summation over the subgroup $\Z[\i]^2$:
$$4\sum_{n\in \Z[\i]}
\ex[(n+a)\tau(n+a)^*+\re((n+a)b^*)]=4\hh_{ab}(2\tau).$$
On the other hand, the summation (\ref{sum-L}) is 
\begin{eqnarray*}
& &\sum_{q\in \F^2}\sum_{n\in \Z[\i]^2}
\begin{array}{l}
\ex[\frac{1}{2}\re({n\over 1+\i}q^*)]
\cdot \ex[({n\over 1+\i}+a)\tau({n\over 1+\i}+a)^*+
\re(({n\over 1+\i}+a)b^*)]
\end{array}
\\
&=&\sum_{q\in \F^2}\ex[-\frac{1}{2}\re(aq^*)]\sum_{n\in \Z[\i]^2}
\begin{array}{l}
\ex[\frac{1}{2}(n+(1+\i)a)\tau(n+(1+\i)a)^*]\\
\times\ex[\re((n+(1+\i)a)({b+q\over 1-\i})^*)]
\end{array}
\\
&=&\sum_{q\in \F^2}\ex[-\frac{1}{2}\re(aq^*)]\hh_{(1+\i)a,{b+q\over 1-\i}}(\tau).
\end{eqnarray*}
For $a\in {1\over 1+\i}\F_2^2$ and $b=(0,0)$, we have the rests.
\qed

\begin{cor}
\label{4term-ex}
\begin{eqnarray*}
\hh_{[0001]}^2(2\tau)+\hh_{[1111]}^2(2\tau)&=&{
\hh_{[0000]}(\tau)+\hh_{[0010]}(\tau)\over 2}
{\hh_{[0001]}(\tau)+\hh_{[0011]}(\tau)\over 2},\\
\hh_{[0010]}^2(2\tau)+\hh_{[1111]}^2(2\tau)&=&
{\hh_{[0000]}(\tau)+\hh_{[0001]}(\tau)\over 2}
{\hh_{[0010]}(\tau)+\hh_{[0011]}(\tau)\over 2},\\
\hh_{[0011]}^2(2\tau)+\hh_{[1111]}^2(2\tau)&=&{
\hh_{[0000]}(\tau)\hh_{[0011]}(\tau)+\hh_{[0001]}(\tau)\hh_{[0010]}(\tau)
\over 2}.
\end{eqnarray*}
\end{cor}
\pr 
By Proposition \ref{hirei} and Pl\"ucker relations, we have 
\begin{eqnarray*}
&&\!\hh_{[0001]}^2(2\tau)\!+\!\hh_{[1111]}^2(2\tau)\!=\!
\hh_{[0000]}^2(2\tau)\!-\!\hh_{[0100]}^2(2\tau),\\
&&\!\hh_{[0010]}^2(2\tau)\!+\!\hh_{[1111]}^2(2\tau)\!=\!
\hh_{[0000]}^2(2\tau)\!-\!\hh_{[1000]}^2(2\tau),\\
&&\!\hh_{[0011]}^2(2\tau)\!+\!\hh_{[1111]}^2(2\tau)\!=\!
\hh_{[0000]}^2(2\tau)\!-\!\hh_{[0100]}^2(2\tau)
\!-\!\hh_{[1000]}^2(2\tau)\!+\!\hh_{[1100]}^2(2\tau).
\end{eqnarray*}
By Theorem \ref{2tau-formula}, we have the corollary.
\qed

\section{Thomae type formula for K3 surfaces}

\subsection{Main Theorem}
\begin{theorem}
\label{main}
Suppose that $x$ is in a neighborhood $U$ of our reference point 
$\dot x\in \Mp(3,6)$.
Let $\tau$ be an element of $\D$ defined in (\ref{periodmatrix}).
Then we have
\begin{equation}
\label{Main-Th}
\hh_{\la J \ra}^2(\tau)=
{1\over 4\pi^4}x\la J \ra \w_{34}(x)^2
\end{equation}
for any $\la J \ra \in P_{3,3}$
\end{theorem}
\begin{remark}
Using the notations in Proposition \ref{series-exp}, 
the above value is equal to $\nu_{34}\la J \ra F_S(z)^2$.
\end{remark}
\pr 
By the first statement of 
Proposition \ref{series-exp},
$x\la J \ra \w_{34}(x)^2$ is a holomorphic function on $D_H$.
We use actions of $\widetilde{\bG}_H(2)$ 
(defined in Proposition \ref{monodromy})
and $U_{22}^{tp}(1+\i)$
on the domains $D_H$ and $\D$ to compare two
functions $\hh_{\la J \ra}^2(\tau)$
and $x\la J \ra \w_{34}(x)^2$. 
\begin{lemma}
Let $\mathcal J_1(g,\tau)$ $(\gamma\in U_{2,2}^{tp}(1+\i),\tau
 \in \D)$
and $\mathcal J_2(R,\w)$ $(R \in \widetilde{\bG}_H(2), \w \in D_H)$
be two cocycles defined by 
\begin{equation}
\label{two cocycles compare}
\mathcal J_1(g,\tau)={\hh_{\la J \ra}^2(g \tau)\over \hh_{\la J
 \ra}^2(\tau)},\quad
\mathcal J_2(R,\w)={R\w_{34}^2 \over \w_{34}^2}.
\end{equation}
Then they coincide via the isomorphisms 
(\ref{fundmental isom on domains})
 and (\ref{isom of discrete group extended}).
\end{lemma}
\pr
Since the group $U^{tp}_{2,2}(1+\i)$ is generated by 
$U_{22}(1+\i)/\la\i E_4\ra$ and $tp$, it is enough to show the identity
(\ref{two cocycles compare}) for 
\begin{itemize}
\item[$(1)$]
$g\in U_{22}(1+\i)/\la \i E_4\ra$ and $R=R_g$, and 
\item[$(2)$]
$g=tp$ and $R=T$.
\end{itemize}
In the case (1), the statement follows from the equality
(\ref{auto-factor}), and that for (2) follows from
Proposition \ref{hirei} (1) and Proposition \ref{monodromy} (4).
\qed

By the above lemma, the function 
$f(x)=\hh_{\la J \ra}^2(\tau(x))/(x\la J \ra
\w_{34}(x)^2)$ becomes a function on $\D/U_{22}^{tp}\simeq X/\la as\ra$.
The space $X/\la as\ra$ can be compactified by the embedding 
$pl^*:X/\la as\ra \to \P^4$.
It is shown in \cite{Ma} that 
the zero of $\hh_{\la J \ra}(\tau(x))^2$ coincides with 
that of $x\la J\ra$. Hence $f(x)$ is a constant map.
We evaluate this constant by taking the degeneration for
$z_2 \to 0$, $z_3\to 0$ in the affine open set of $\overline{X}$
defined by
$$
\biggl\{
\pmatrix{
1 & 1 & 0 & 0 & 1 & 1 \cr
-z_1&-1& 1 & 0 & 0 &-z_3\cr
-z_2&0 & 0 & 1 &-1 &-z_4\cr}\mid   z_1, \dots, z_4 \in \C
\biggr\}.
$$
Under this limit, we have
$\w_{13}(x), \w_{24}(x) \to 0$ and
$$
\w_{34} \to 
2\w_B(z_1)\w_B(z_4),\quad
\w_{14} \to 
2\w_A(z_1)\w_B(z_4),\quad
\w_{23} \to 
2\w_B(z_1)\w_A(z_4),
$$
where $\w_A(\lambda)$, $\w_B(\lambda)$ are defined in 
(\ref{eq:definition of elliptic period}).
We set $\w_{ij}=\lim_{z_2,z_3\to 0}\w_{ij}(x)$.
Then we have
$$
\lim_{z_2,z_3 \to 0}\tau(x)= 
\diag ({\w_{14}\over \w_{34}},{\w_{23}\over \w_{34}}),
$$
and by $x\la 135 \ra\to 1$,
\begin{eqnarray*}
\lim_{z_2,z_3 \to 0} f(x)&=& \hh_{[0000]}^2
(\diag ({\w_{14}\over \w_{34}},{\w_{23}\over \w_{34}})) 
/(x\la 135 \ra \w_{34}^2)
\\
&=&
\vartheta_{[00]}^4({\w_{A}(z_1)\over \w_{B}(z_1)})
\vartheta_{[00]}^4({\w_{A}(z_4)\over \w_{B}(z_4)})
/(4\w_{B}^2(z_1)\w_{B}^2(z_4)) \\
&=&
{1 \over 4\pi^4}
\end{eqnarray*}
by Jacobi's formula (\ref{Jacobi's formula}).
\qed

\section{Mean iterations}

\subsection{Mean iteration associated to $D_4$ degeneration}
In this and next subsections, 
we apply the main identity (\ref{Main-Th})
to the study
of mean iterations. In this subsection, we consider 
configurations of six lines which contain
three lines intersecting at one point. In this degeneration,
three $A_1$ singularities on $\widehat{\cS}^*$ confluent to one $D_4$ singularity.
This degeneration is obtained by taking the limit $x\la J\ra \to 0$.
We consider the case $\la J\ra =\la 123\ra$.

\begin{proposition}
\label{x123=0}
The two preimages of the map $\overline{pl}$ on the subvariety defined by
$x\la123\ra=0$ are expressed as
$$\pmatrix{1 & 1& 0&  0& 1& 1\cr
-{x\la124\ra-x\la125\ra-x\la134\ra+x\la135\ra\over -x\la134\ra+x\la135\ra}
& -1& 1& 0& 0
& -{x\la124\ra x\la135\ra-x\la125\ra x\la134\ra\over 
(-x\la134\ra+x\la135\ra)x\la125\ra}\cr 
0& 0& 0& 1& -1& -{x\la125\ra-x\la124\ra\over x\la125\ra}\cr},
$$
$$\pmatrix{1 & 1& 0&  0& 1& 1\cr
-{x\la134\ra-x\la124\ra\over x\la134\ra}& -1& 1& 0& 0& 0\cr 
-{x\la135\ra x\la124\ra- x\la125\ra x\la134\ra \over 
(x\la135\ra-x\la125\ra) x\la134\ra }& 0& 0& 1& -1& 
-{x\la124\ra-x\la125\ra-x\la134\ra+x\la135\ra\over 
x\la135\ra-x\la125\ra}\cr}.
$$
\end{proposition}

\label{our AGM}
Let $m$ be a map from $(\Rp)^4$ to $(\Rp)^4$ given by 
\begin{equation}
\label{def:new iteration}
m:(\Rp)^4\ni u=(u_1,\dots,u_4)\mapsto 
(m_1(u),\dots,m_4(u))\in (\Rp)^4,
\end{equation}
where 
\begin{eqnarray*}
m_1(u)&=&{u_1+u_2+u_3+u_4\over 4},\hspace{10mm}
m_2(u)={\sqrt{(u_1+u_3)(u_2+u_4)}\over 2},\\
m_3(u)&=&{\sqrt{(u_1+u_2)(u_3+u_4)}\over 2},\quad
m_4(u)=\sqrt{u_1u_4+u_2u_3\over 2}.
\end{eqnarray*}
For an element $c=(c_1,\dots,c_4)\in (\Rp)^4$ with $c_1>c_2>c_3>c_4$,  
we define a vector valued sequence $\{m^n(c)=(m_1^n(c),\dots, m_4^n(c))\}_{n\in \N}$
by
$$
m^n(c)=\overbrace{m\circ\cdots\circ m}^{n}(c).
$$
\begin{lemma}
\label{converge}
\begin{itemize}
\item[$(1)$]
The components of the sequence $\{(m^n_1(c),\dots, m^n_4(c))\}_{n\in \N}$ 
converge and 
have a common limit $m^\infty_*(c)$. The convergence is quadratic. 
\item[$(2)$]
$$\lim_{n\to \infty}
{m^n_1(c)^2-m^n_2(c)^2\over m^n_3(c)^2-m^n_4(c)^2}
=\lim_{n\to \infty}
{m^n_1(c)^2-m^n_3(c)^2\over m^n_2(c)^2-m^n_4(c)^2}=1.
$$
\end{itemize}
\end{lemma}
\pr (1)  Since 
\begin{eqnarray*}
& &m_1(c)^2-m_2(c)^2={(c_1-c_2+c_3-c_4)^2\over 16},\\
& &m_2(c)^2-m_3(c)^2={(c_1-c_4)(c_2-c_3)\over 4},\\
& &m_3(c)^2-m_4(c)^2={(c_1-c_2)(c_3-c_4)\over 4},
\end{eqnarray*}
we have 
$m_1(c)>m_2(c)>m_3(c)>m_4(c)$ for $c_1>c_2>c_3>c_4>0.$
We can easily see that 
$$c_4<m_4(c)<\cdots<m_4^n(c)<m_1^n(c)<\cdots<m_1(c)<c_1,$$ 
the sequences $\{m_1^n(c)\}$ and $\{m_4^n(c)\}$ converge.
We set $\mu_1=\lim\limits_{n\to\infty}m_1^n(c)$ and 
$\mu_4=\lim\limits_{n\to\infty}m_4^n(c)$. 
Since 
\begin{eqnarray}
\nonumber
m_1(c)^2-m_4(c)^2&=&m_1(c)^2-m_2(c)^2+m_2(c)^2-m_4^2(c)\\
\nonumber
&=&{(c_1-c_2+c_3-c_4)^2\over 16}+{(c_1-c_3)(c_2-c_4)\over 4}\\
\label{inequality for quadratic converge}
&<&{(c_1-c_4)^2\over 4}+{(c_1-c_4)^2\over 4}={(c_1-c_4)^2\over 2},
\end{eqnarray}
we have $\mu_1^2-\mu_4^2\leq\frac{1}{2}(\mu_1-\mu_4)^2$.
If $\mu_1>\mu_4$, then $\mu_1+\mu_4\leq \frac{1}{2}(\mu_1-\mu_4)$,
which implies $\mu_1+3\mu_4\leq 0$. This is a contradiction.

By the inequality (\ref{inequality for quadratic converge}),
the convergence of $m_1(c)-m_4(c)={m_1(c)^2-m_4(c)^2\over m_1(c)+m_4(c)}$,
is quadratic.

\smallskip\noindent (2)  
We have 
$${m_1(c)^2-m_2(c)^2\over m_3(c)^2-m_4(c)^2}={1\over 4}\left(
{c_1-c_2 \over c_3-c_4}+2+{c_3-c_4\over c_1-c_2}
\right)\ge1.$$ 
We set 
$$s_n={m_1^n(c)^2-m_2^n(c)^2\over m_3^n(c)^2-m_4^n(c)^2},\quad 
t_n={m_1^n(c)+m_2^n(c)\over m_3^n(c)+m_4^n(c)},\quad 
f(s,t)={1\over 4}({s\over t}+2+{t\over s}).$$
Then $s_n$ and $t_n$ satisfy 
$$s_{n+1}=f(s_n,t_n),\quad s_n,t_n\ge 1,\quad \lim_{n\to \infty}t_n=1.$$
Note that $f(hs,ht)=f(s,t)$ for any $h\in \Rp$, $f(s,s)=1$ and that
$$f(s,1)=f(1,s)={1\over 4}(s+2+{1\over s})<{1\over 4}(s+2s+ s)<s$$
for any $s>1$.
If $s_n>t_n$ then $(s_n/t_n)>1$ and 
$$s_{n+1}=f(s_n,t_n)=f({s_n\over t_n},1)<{s_n\over t_n}<s_n.$$
If $s_n\le t_n$ then $(t_n/s_n)\ge 1$ and 
$$s_{n+1}=f(s_n,t_n)=f({t_n\over s_n},1)\le {t_n\over s_n}\le t_n.$$
Thus we have $s_{n+1}\le \max(s_n,t_n)$. Since 
$\lim\limits_{n\to\infty}t_n=1$, for any $\e>0$ there exists 
$N\in \N$ such that $t_n<1+\e$ for any $n>N$. If there exists $n_0>N$
such that $s_{n_0}\le t_{n_0}$, then 
$s_n<1+\e$ for any $n\ge n_0$; this means $\lim\limits_{n\to\infty}s_n=1$.
Otherwise, i.e. $s_n>t_n$ for any $n>N$, then 
$s_n$ is monotonously decreasing. Thus the limit
$\lim\limits_{n\to\infty}s_n$ exists.
Let $n\to \infty$ for $s_{n+1}=f(s_n,t_n)$, then 
we have  $\lim\limits_{n\to\infty}s_n=1.$
Similarly we can show 
$\lim\limits_{n\to \infty}
{m^n_1(c)^2-m^n_3(c)^2\over m^n_2(c)^2-m^n_4(c)^2}=1
$. 
\qed

\begin{theorem}
\label{AGM4}
The common limit $m^\infty_*(c)$ can be expressed 
as 
$$
m^\infty_*(c)=
\sqrt{c_1^2-c_2^2\over c_3^2-c_4^2}
{c_3\over F_S(z)} 
=\sqrt{c_1^2-c_3^2\over c_2^2-c_4^2}{c_2\over F_S(w)},
$$
where $z=\pmatrix{z_1 & z_3\cr z_2 & z_4}$ and 
$w=\pmatrix{w_1 & w_3\cr w_2 & w_4}$ are given as
$$
z=\pmatrix{1-{c_3^2-c_4^2\over c_1^2-c_2^2}& 
1-{c_1^2(c_3^2-c_4^2)\over c_3^2(c_1^2-c_2^2)}\cr 
0 & 1-{c_4^2\over c_3^2}},\quad 
w=\pmatrix{1-{c_4^2\over c_2^2}& 0 \cr
1-{c_1^2(c_2^2-c_4^2)\over c_2^2(c_1^2-c_3^2)}&
1-{c_2^2-c_4^2\over c_1^2-c_3^2}}.
$$

\end{theorem}
\begin{remark}
For a given $c=(c_1,\dots,c_4)$, the hypergeometric series $F_S$ 
in Theorem \ref{AGM4} may not converge. 
By Lemma \ref{converge} (2), there exists $n\in \N$ 
such that it converges for $m^n(c)$ instead of $c$.
\end{remark}
\pr 
There exists $\tau\in \D$  such that $\hh_{[1111]}(\tau)=0$ and 
$$\hh_{[0000]}(\tau):\hh_{[0001]}(\tau):\hh_{[0010]}(\tau):\hh_{[0011]}(\tau)
=c_1:c_2:c_3:c_4.$$
By Corollary \ref{4term-ex}, we have 
\begin{eqnarray*}
& &(\hh_{[0000]}(2\tau),\hh_{[0001]}(2\tau),
\hh_{[0010]}(2\tau),\hh_{[0011]}(2\tau))\\
&=&
m(\hh_{[0000]}(\tau),\hh_{[0001]}(\tau),
\hh_{[0010]}(\tau),\hh_{[0011]}(\tau)),
\end{eqnarray*}
since $\hh_{[1111]}(\tau)=0$.
By the homogeneity of $m_1,\dots,m_4$, $m^\infty_*$ satisfies 
$$m^\infty_*(c)=c_1m^\infty_*(1,{c_2\over c_1},{c_3\over c_1},{c_4\over c_1}).
$$
Thus 
\begin{eqnarray*}
& &m^\infty_*(c)=c_1m^\infty_*(1,{\hh_{[0001]}(\tau)\over \hh_{[0000]}(\tau)},
{\hh_{[0010]}(\tau)\over \hh_{[0000]}(\tau)},
{\hh_{[0011]}(\tau)\over \hh_{[0000]}(\tau)})\\
&=&{c_1\over\hh_{[0000]}(\tau)} 
m^\infty_*(\hh_{[0000]}(\tau),\hh_{[0001]}(\tau),
\hh_{[0010]}(\tau),\hh_{[0011]}(\tau))\\
&=&{c_1\over\hh_{[0000]}(\tau)} 
m^\infty_*(m(\hh_{[0000]}(\tau),\hh_{[0001]}(\tau),
\hh_{[0010]}(\tau),\hh_{[0011]}(\tau)))\\
&=&{c_1\over\hh_{[0000]}(\tau)} 
m^\infty_*(\hh_{[0000]}(2\tau),\hh_{[0001]}(2\tau),
\hh_{[0010]}(2\tau),\hh_{[0011]}(2\tau))\\
&=&{c_1\over\hh_{[0000]}(\tau)} 
m^\infty_*(\hh_{[0000]}(2^n\tau),\hh_{[0001]}(2^n\tau),
\hh_{[0010]}(2^n\tau),\hh_{[0011]}(2^n\tau))\\
&\to&{c_1\over\hh_{[0000]}(\tau)} \quad \textrm{as}\quad n\to \infty,
\end{eqnarray*}
since $\hh_{[ab]}(2^n\tau)$ converge to $1$ for $a=(0,0)$ and 
any $b\in \F_2^2$.  
By Proposition \ref{x123=0}, the preimages of 
$$[x\la 123\ra,x\la 124\ra,x\la 125\ra,x\la 134\ra,x\la 135\ra]
=[0,c_4^2,c_3^2,c_2^2,c_1^2]$$
for the map $\overline{pl}:\overline{X}\to  \P^4\subset \P^9$ are given by 
$$x=\pmatrix{
1 & 1 & 0 & 0 & 1 & 1\cr
-z_1& -1&1& 0 & 0 &-z_3\cr
0& 0&1& -1 & 0 &-z_4\cr}, $$
$$z_1=1-{c_3^2-c_4^2\over c_1^2-c_2^2}, \quad z_2=0,\quad
z_3=1-{c_1^2(c_3^2-c_4^2)\over c_3^2(c_1^2-c_2^2)}, \quad
z_4=1-{c_4^2\over c_3^2},
$$
with
$$\sqrt{x\la 135\ra}={c_1\over c_3}\sqrt{c_3^3-c_4^2\over c_1^2-c_2^2},$$
and 
$$
x=\pmatrix{
1 & 1 & 0 & 0 & 1 & 1\cr
-z_1& -1&1& 0 & 0 & 0\cr
-z_2& 0&1& -1 & 0 &-z_4\cr}, $$
$$z_1=1-{c_4^2\over c_2^2},\quad 
z_2=1-{c_1^2(c_2^2-c_4^2)\over c_2^2(c_1^2-c_3^2)},\quad z_3=0,\quad 
z_4=1-{c_2^2-c_4^2\over c_1^2-c_3^2},
$$
with 
$$\sqrt{x\la 135\ra}={c_1\over c_2}\sqrt{c_2^2-c_4^2\over c_1^2-c_3^2}.$$
Theorem \ref{main} implies this theorem. 
\qed
\subsection{Mean iteration associated to Kummer locus}
\label{Borchardt's AGM}
The Borchardt's mean iteration is obtained from
the restriction of Thomae type formula for K3 surfaces to
the Kummer locus.
In this subsection, we explain how to recover limit formulas
in \cite{B}, \cite{MT} and \cite{Me}
from our main theorem. 
\begin{proposition}
\label{Igusa-locus}
Let $c_1>c_2>c_3>c_4$ be real numbers such that
$c_1-c_2-c_3+c_4 > 0$.
We set
\begin{eqnarray*}
Q_1&=&
(c_1\!+\! c_2\!+\! c_3\!+\! c_4)(c_1\!+\! c_2\!-\! c_3\!-\! c_4)
(c_1\!-\! c_2\!+\! c_3\!-\! c_4)(c_1\!-\! c_2\!-\! c_3\!+\! c_4),\\
\end{eqnarray*}
and
\begin{eqnarray}
\label{kummer locus coordinate z}
\pmatrix{z_1 & z_3 \cr z_2 &z_4}=\pmatrix{
1-{c_4(c_1^2-c_2^2+c_3^2-c_4^2-\sqrt{Q_1})\over 2c_2(c_1c_3-c_2c_4)}& 
1-{c_1(c_1^2-c_2^2+c_3^2-c_4^2-\sqrt{Q_1})\over 2c_3(c_1c_3-c_2c_4)}\cr 
1-{c_1(c_1^2+c_2^2-c_3^2-c_4^2-\sqrt{Q_1})\over 2c_2(c_1c_2-c_3c_4)}& 
1-{c_4(c_1^2+c_2^2-c_3^2-c_4^2-\sqrt{Q_1})\over 2c_3(c_1c_2-c_3c_4)}}.
\end{eqnarray}
Then
\begin{equation}
\label{definition of x in kummer}
x=\pmatrix{
1 & 1 & 0 & 0 & 1 & 1 \cr
-z_1&-1& 1 & 0 & 0 &-z_3\cr
-z_2&0 & 0 & 1 &-1 &-z_4\cr}
\end{equation}
lies on the Kummer locus.
In this case, we have
$$
[x\la123\ra: 
x\la135\ra:
x\la134\ra:
x\la125\ra: 
x\la124\ra]
=
[c_0^2:
c_1^2:
c_2^2:
c_3^2:
c_4^2],
$$
where
$$
c_0^2={c_1^2-c_2^2-c_3^2+c_4^2+\sqrt{Q_1}\over 2}.
$$

\end{proposition}

Let $m$ be a map from $(\Rp)^4$ to $(\Rp)^4$ given by 
\begin{equation}
\label{def:kummer iteration map}
m:(\Rp)^4\ni u=(u_1,\dots,u_4)\mapsto 
(m_1(u),\dots,m_4(u))\in (\Rp)^4,
\end{equation}
where 
\begin{eqnarray*}
m_1(u)&=&{u_1+u_2+u_3+u_4\over 4},\quad
m_2(u)={\sqrt{u_1u_2}+\sqrt{u_3u_4}\over 2},\\
m_3(u)&=&{\sqrt{u_1u_3}+\sqrt{u_2u_4}\over 2},\hspace{8mm}
m_4(u)={\sqrt{u_1u_4}+\sqrt{u_2u_3}\over 2}.
\end{eqnarray*}
Note that if $u_1>u_2>u_3>u_4$ then 
$$m_1(u)>m_2(u)>m_3(u)>m_4(u).$$
For an element $c=(c_1,\dots,c_4)\in (\Rp)^4$ with $c_1>c_2>c_3>c_4$,  
we define a vector valued sequence $\{m^n(c)=(m_1^n(c),\dots, m_4^n(c))\}_{n\in \N}$ by
$$
m^n(c)=\overbrace{m\circ\cdots\circ m}^{n}(c).
$$

In \cite{B} and \cite{Me}, they 
prove that
the common limit $m^\infty_*(c)$ is expressed in terms of period integrals of 
a hyperelliptic curve of genus 2. In \cite{MT}, they give its expression
in terms of
the period integral $\w_{34}(x)$ of the \KK surface $\cS(x)$. 
Here, we give its expression by the hypergeometric series $F_S$.

\begin{theorem}
\label{BAGM4}
We can express the common limit $m^\infty_*(c)$ by
\begin{eqnarray*}
m^\infty_*(c)
&=&{4\sqrt{c_2c_3(c_1c_2-c_3c_4)(c_1c_3-c_2c_4)}
\over (\sqrt{d_1d_2}-\sqrt{d_3d_4})(\sqrt{d_1d_3}-\sqrt{d_2d_4})}
{1\over F_S(z)}\\
&=&{4\sqrt{c_2c_3(c_1c_2-c_3c_4)(c_1c_3-c_2c_4)}
\over (\sqrt{d_1d_2}+\sqrt{d_3d_4})(\sqrt{d_1d_3}+\sqrt{d_2d_4})}
{1\over F_S(w)}
\end{eqnarray*}
where 
$$z=\pmatrix{z_1 & z_3 \cr z_2 &z_4}=\pmatrix{
1-{c_4\over c_2}{\sqrt{d_1d_3}-\sqrt{d_2d_4}\over\sqrt{d_1d_3}+\sqrt{d_2d_4}}&
1-{c_1\over c_3}{\sqrt{d_1d_3}-\sqrt{d_2d_4}\over\sqrt{d_1d_3}+\sqrt{d_2d_4}}
\cr
1-{c_1\over c_2}{\sqrt{d_1d_2}-\sqrt{d_3d_4}\over\sqrt{d_1d_2}+\sqrt{d_3d_4}}&
1-{c_4\over c_3}{\sqrt{d_1d_2}-\sqrt{d_3d_4}\over\sqrt{d_1d_2}+\sqrt{d_3d_4}}
\cr},
$$
$$w=\pmatrix{z_1 & z_3 \cr z_2 &z_4}=\pmatrix{
1-{c_4\over c_2}{\sqrt{d_1d_3}+\sqrt{d_2d_4}\over\sqrt{d_1d_3}-\sqrt{d_2d_4}}&
1-{c_1\over c_3}{\sqrt{d_1d_3}+\sqrt{d_2d_4}\over\sqrt{d_1d_3}-\sqrt{d_2d_4}}
\cr
1-{c_1\over c_2}{\sqrt{d_1d_2}+\sqrt{d_3d_4}\over\sqrt{d_1d_2}-\sqrt{d_3d_4}}&
1-{c_4\over c_3}{\sqrt{d_1d_2}+\sqrt{d_3d_4}\over\sqrt{d_1d_2}-\sqrt{d_3d_4}}
\cr},
$$
$$d_1=c_1+c_2+c_3+c_4,\quad d_2=c_1+c_2-c_3-c_4,$$
$$d_3=c_1-c_2+c_3-c_4,\quad d_4=c_1-c_2-c_3+c_4.$$
\end{theorem}
\pr Let $c_1,\dots,c_4$ be elements in $\Rp$ with $c_1>\cdots>c_4$.
Though the value $c_1-c_2-c_3+c_4$ may be negative, by the inequality
$$
m_1(c)\!-\! m_2(c)\!-\! m_3(c) \! +\! m_4(c)=
{1\over 4}(\sqrt{c_1}\!-\! \sqrt{c_2}\!-\! \sqrt{c_3} \! +\!
\sqrt{c_4})^2>0,
$$
we may assume $c_1-c_2-c_3+c_4\ge 0$ by applying the map $m$.
Let 
$\left(
\begin{array}{cc}z_1,z_2 \\z_3,z_4
\end{array}\right)
$ be $2\times 2$ matrix defined as 
(\ref{kummer locus coordinate z}).
Then $x$ defined in (\ref{definition of x in kummer})
 lies on the Kummer locus.
In this case, the theta constants $\hh_{[a,b]}(\tau)$
coincide with the square of Riemann's theta constants.
Using $2\tau$-formulas for Riemann's theta constants, and similar
argument as in Theorem \ref{AGM4},
we have 
\begin{eqnarray*}
m_*^\infty(c)
={c_1\over \hh_{[0000]}(\tau)}.
\end{eqnarray*}
By Theorem \ref{main} and Proposition \ref{series-exp},
we have the first expression of $m_*^\infty(c)$ 
by $F_S$.
By putting 
$$c_0^2={c_1^2-c_2^2-c_3^2+c_4^2-\sqrt{Q_1}\over 2},$$
we have the other expression of  $m_*^\infty(c)$.  
\qed

\section{Functional equations for $F_S$}
\label{FE}
The common limit $m_*^\infty(c)$ in Theorem \ref{AGM4} (resp. \ref{BAGM4}) 
satisfies 
$$ m_*^\infty(c)=m_*^\infty(m(c)).$$
This property implies functional equations for 
the hypergeometric function $F_S$. 
\begin{theorem}
\label{FE1}
We have the following functional equations for $F_S$:
\comment{
\begin{eqnarray*}
F_S\left({(c_1\!-\! c_2\!-\! c_3\!+\! c_4)^2\over 
(c_1\!-\! c_2\!+\! c_3\!-\! c_4)^2},0,
{2(c_1c_4\!-\! c_2c_3)(c_1^2\!-\! c_2^2\!-\! c_3^2\!+\! c_4^2)\over 
(c_1\!+\! c_2)(c_3\!+\! c_4)(c_1\!-\! c_2\!+\! c_3\!-\! c_4)^2},
{(c_1\!-\! c_2)(c_3\!-\! c_4)\over (c_1\!+\! c_2)(c_3\!+\! c_4)}\right)\\
={(c_1\!-\! c_2\!+\! c_3\!-\! c_4)(c_3\!+\! c_4)\over 4(c_1\!-\! c_2)c_3}
F_S\left({c_1^2\!-\! c_2^2\!-\! c_3^2\!+\! c_4^2\over c_1^2\!-\! c_2^2},0,
{c_1^2c_4^2\!-\! c_2^2c_3^2\over (c_1^2\!-\! c_2^2)c_3^2},
{c_3^2\!-\! c_4^2\over c_3^2}\right).
\end{eqnarray*}
}
\begin{eqnarray*}
F_S(m(z))&=&
{(c_1\!-\! c_2\!+\! c_3\!-\! c_4)(c_3\!+\! c_4)\over 4(c_1\!-\! c_2)c_3}
F_S(z),\\
F_S(m(w))&=&
{(c_1\!+\! c_2\!-\! c_3\!-\! c_4)(c_2\!+\! c_4)\over 4(c_1\!-\! c_3)c_2}
F_S(w),\end{eqnarray*}
where $z$ and $w$ are given in Theorem \ref{AGM4} and 
\begin{eqnarray*}
m(z)&=&\pmatrix{
{(c_1\!-\! c_2\!-\! c_3\!+\! c_4)^2\over (c_1\!-\! c_2\!+\! c_3\!-\! c_4)^2}
&
{2(c_1c_4\!-\! c_2c_3)(c_1^2\!-\! c_2^2\!-\! c_3^2\!+\! c_4^2)\over 
(c_1\!+\! c_2)(c_3\!+\! c_4)(c_1\!-\! c_2\!+\! c_3\!-\! c_4)^2}
\cr
0
&
{(c_1\!-\! c_2)(c_3\!-\! c_4)\over (c_1\!+\! c_2)(c_3\!+\! c_4)}},\\
m(w)&=&\pmatrix{
{(c_1\!-\! c_3)(c_2\!-\! c_4)\over (c_1\!+\! c_3)(c_2\!+\! c_4)}
&
0
\cr
{2(c_1c_4\!-\! c_2c_3)(c_1^2\!-\! c_2^2\!-\! c_3^2\!+\! c_4^2)\over 
(c_1\!+\! c_3)(c_2\!+\! c_4)(c_1\!+\! c_2\!-\! c_3\!-\! c_4)^2}
&
{(c_1\!-\! c_2\!-\! c_3\!+\! c_4)^2\over (c_1\!+\! c_2\!-\! c_3\!-\! c_4)^2}}.
\end{eqnarray*}
\end{theorem}

\begin{theorem}
\label{FE2}
We have the following functional equations for $F_S$:
\comment{
\begin{eqnarray*}
& & {4\sqrt{(\sqrt{c_1c_2}-\sqrt{c_3c_4})(\sqrt{c_1c_3}-\sqrt{c_2c_4})}\over
\sqrt{(\sqrt{c_1c_2}+\sqrt{c_3c_4})(\sqrt{c_1c_3}+\sqrt{c_2c_4})
d_2d_3}} F_S(m(z))\\
&=&{(\sqrt{d_1d_2}-\sqrt{d_3d_4})(\sqrt{d_1d_3}-\sqrt{d_2d_4})
\over 4\sqrt{c_2c_3(c_1c_2-c_3c_4)(c_1c_3-c_2c_4)}}
 F_S(z),
\end{eqnarray*}
}
\begin{eqnarray*}
F_S(m(z))&=&{1\over 16}
{\sqrt{d_2d_3}
(\sqrt{d_1d_2}-\sqrt{d_3d_4})(\sqrt{d_1d_3}-\sqrt{d_2d_4})
\over 
\sqrt{c_2c_3}(\sqrt{c_1c_2}-\sqrt{c_3c_4})(\sqrt{c_1c_3}-\sqrt{c_2c_4})}
 F_S(z),
\\
F_S(m(w))&=&{(\sqrt{d_1d_2}+\sqrt{d_3d_4})(\sqrt{d_1d_3}+\sqrt{d_2d_4})
\over 4\sqrt{c_2c_3d_2d_3}}
 F_S(w),
\end{eqnarray*}
where $z$ and $w$ are given in Theorem \ref{BAGM4} and 
\begin{eqnarray*}
m(z)&=&\pmatrix{
1-{2(\sqrt{c_1c_4}+\sqrt{c_2c_3})(\sqrt{c_1c_3}-\sqrt{c_2c_4})
\over(\sqrt{c_1c_2}+\sqrt{c_3c_4})d_3}&
1-{(\sqrt{c_1c_3}-\sqrt{c_2c_4})d_1
\over(\sqrt{c_1c_3}+\sqrt{c_2c_4})d_3}
\cr
1-{(\sqrt{c_1c_2}-\sqrt{c_3c_4})d_1
\over(\sqrt{c_1c_2}+\sqrt{c_3c_4})d_2}
&
1-{2(\sqrt{c_1c_4}+\sqrt{c_2c_3})(\sqrt{c_1c_2}-\sqrt{c_3c_4})
\over(\sqrt{c_1c_3}+\sqrt{c_2c_4})d_2}
\cr},
\\
m(w)&=&\pmatrix{
1-{(\sqrt{c_1c_4}+\sqrt{c_2c_3})d_3
\over2(\sqrt{c_1c_2}+\sqrt{c_3c_4})(\sqrt{c_1c_3}-\sqrt{c_2c_4})}&
1-{d_1d_3\over4(c_1c_3-c_2c_4)}
\cr
1-{d_1d_2\over4(c_1c_2-c_3c_4)}
&
1-{(\sqrt{c_1c_4}+\sqrt{c_2c_3})d_2
\over2(\sqrt{c_1c_3}+\sqrt{c_2c_4})(\sqrt{c_1c_2}-\sqrt{c_3c_4})}
\cr}.
\end{eqnarray*}
\end{theorem}
\pr
To obtain the expression of $m(z)$ and $m(w)$, we use the equalities:
$$
m_1(c)  +\epsilon_1  m_2(c)  + \epsilon_2  m_3(c) 
+\epsilon_1\epsilon_2  m_4(c)=
{1\over 4}(\sqrt{c_1} +\epsilon_1 \sqrt{c_2} +\epsilon_2 \sqrt{c_3} 
+\epsilon_1\epsilon_2 \sqrt{c_4})^2
$$
for $\epsilon_1, \epsilon_2=\pm 1$ and
\begin{eqnarray*}
& &m_i(c)m_j(c)\!-\! m_k(c)m_l(c)=
{1\over 8}(\sqrt{c_ic_j}\!-\!\sqrt{c_kc_l})(c_i \!+\! c_j\!-\! c_k\!-\!
c_l)
\end{eqnarray*}
for $(i,j,k,l)=(1,2,3,4), (1,3,2,4), (1,4,2,3)$.
\qed

\bigskip
\begin{flushleft}
\begin{minipage}{5.8cm}
Keiji \textsc{Matsumoto}\\
Department of Mathematics\\
Hokkaido University\\
Sapporo 060-0810, Japan\\
e-mail:matsu@math.sci.hokudai.ac.jp
\end{minipage}
\end{flushleft}

\begin{flushleft}
\begin{minipage}{8.0cm}
Tomohide \textsc{Terasoma}\\
Graduate School of Mathematical Sciences\\
The University of Tokyo \\
Komaba, Meguro, Tokyo, 153-8914, Japan\\
e-mail:terasoma@ms.u-tokyo.ac.jp
\end{minipage}
\end{flushleft}


\begin{thebibliography}{99}
%
\bibitem[B]{B} 
C.W. Borchardt, 
\"Uber das arithmetisch-geometrische Mittel aus vier Elementen, 
\textsl{ Berl. Monatsber,} 53 (1876), 611-621.
%
\bibitem[F]{F} E. Freitag, 
Modulformen zweiten Grades zum rationalen und Gau\ss schen Zahlk\" orper, 
\textsl{Sitzungsber. Heidelb. Akad. Wiss.}, \textbf{1} (1967), 1--49.
%
\bibitem[I]{I}
J. Igusa, 
\textsl{Theta functions,}
Die Grundlehren der mathematischen Wissenshaften in 
Einzeldarstellungen 194,
Springer-Berlin-Heidelberg, New York, 1972.
%
\bibitem[IKSY]{IKSY} K. Iwasaki, H. Kimura, S. Shimomura and M. Yoshida,
\textsl{From Gauss to Painlev\'e},
Vieweg, Braunschweig, Wiesbaden, 1991.
%
\bibitem[KY]{KY} M. Kita and M. Yoshida, 
Intersection theory for twisted cycles I, 
\textsl{Math. Nachr.} {\bf 166} (1994), 287--304.
%
\bibitem[KaM]{KaM} T. Kato and K. Matsumoto,
The common limit of a quadruple sequence and 
the hypergeometric function $F_D$ of three variables,
\textsl{Nagoya Math. J.} {\bf 195} (2009). 113--124.
%
\bibitem[KiM]{KiM}
M. Kita and K. Matsumoto,  
Duality for hypergeometric functions and invariant Gauss-Manin systems, 
\textsl{Compositio Math.}  \textbf{108} (1997), 77--106.
%
%
%
%
\bibitem[MSY]{MSY} K. Matsumoto, T. Sasaki and M. Yoshida, 
The monodromy of the period map of a 4-parameter family of \KK 
surfaces and the Aomoto-Gel'fand hypergeometric function of type (3,6), 
\textsl{Internat. J. of Math.}, \textbf{3} (1992), 1--164.
%
\bibitem[MT]{MT} K. Matsumoto and T. Terasoma, 
Arithmetic-geometric means for hyperelliptic curves and 
Calabi-Yau varieties, to appear in \textsl{Internat. J. of Math.}
%
\bibitem[MY]{MY} K. Matsumoto and M. Yoshida, 
Invariants for some real hyperbolic groups, 
\textsl{Internat. J. of Math.}, \textbf{13} (2002), 415--443.
%
\bibitem[Ma]{Ma} K. Matsumoto,
Theta functions on the bounded symmetric domain of type $I_{2,2}$ and the 
period map of $4$-parameter family of K3 surfaces, 
\textsl{Math. Ann.}, \textbf{295} (1993), 383--408.
%
\bibitem[Me]{Me}
J. Mestre, 
Moyenne de Borchardt et integrales elliptiques,
\textit{C. R. Acad. Sci. Paris Ser. I Math.} 313 (1991), no. 5, 273--276. 
%
\bibitem[Mu]{Mu}
D. Mumford, 
\textit{Tata lectures on Theta I,}
progress in Math 28. Birkh\"auser, Boston-Basel-Berlin, 1983.
%
\bibitem[Te]{T}
T. Terasoma,
 Exponential Kummer coverings and determinants of hypergeometric functions. 
\textsl{Tokyo J. Math.} 16 (1993), no. 2, 497--508.  
%
\bibitem[To]{To}
J. Thomae,
Beitrag zur Bestimmung von $\theta(0, 0, ... , 0)$
durch die Klassenmoduln algebraischer Funktionen,
{\sl  J. Reine Angew. Math.} \textbf{71} (1870), 201--222.
%
\bibitem[Y]{Y} M. Yoshida, \textsl{Hypergeometric Functions, My Love,} 
Aspects of Mathematics, E32, Friedr Vieweg \& Sohn, Braunschweig, 1997. 

\end{thebibliography}
\end{document}